\documentclass[12pt]{article}
\usepackage[english]{babel}
\usepackage{amssymb,latexsym,amsmath,amscd}

\begin{document}

\title
{About isolators of finitely generated subgroups of free groups}
\author
{D.\, Moldavanskii}
\date{}
\maketitle

\begin{abstract}
It is known that in any free group the isolator of finitely generated subgroup is finitely generated subgroup. 
A very simple proof of this statement is proposed.
\end{abstract}

Recall that a subgroup $H$ of some group $G$ is said to be {\it isolated} if for any element $g$ of $G$ and for any positive integer $m$ the inclusion $g^{m}\in H$ implies that $g\in H$. The intersection of all isolated subgroups of group $G$ containing $H$ is called the {\it isolator of subgroup $H$} and is denoted by the symbol $I(H)$. In other words, $I(H)$ is the smallest of isolated subgroups of the group $G$ containing $H$

These notions were introduced by P.\,G.\,Kontorovich \cite{Kont} and then generalized in the work of V.\,N.\,Bezverkhnii \cite{Bezv} as follows:

A subgroup $H$ of some group $G$ is said to be {\it $m$-isolated} for some positive integer $m$ if for any element $g\in G$ the inclusion $g^{m}\in H$ implies that $g\in H$. Let $\pi$ be a nonempty set of primes. Subgroup $H$ of a group $G$ is said to be {\it $\pi$-isolated} if it is $m$-isolated for every $\pi$-number $m$ (or if, equivalently, $H$ is $p$-isolated for every prime $p\in\pi$). The smallest of the $\pi$-isolated subgroups of the group $G$ containing $H$ is called the {\it $\pi$-isolator of subgroup $H$} and is denoted by $I_{\pi}(H)$.

(It should be noted that the original definition of $\pi$-isolation of subgroups formulated in \cite{Bezv} means, in accordance with the definition given here (and consistent with the terminology adopted in modern publications), $\pi'$\nobreakdash-isolation, where $\pi'$ is the complement of set $\pi $ in the set of all primes.)

The following statement is true:\medskip

{\bf Theorem.} {\it For any free group $F$ and for arbitrary nonempty set $\pi$ of primes, the $\pi$-isolator $I_{\pi} (H) $ of every finitely generated subgroup $H$ of  $F$ is a finitely generated subgroup.}\medskip

In the case when the set $\pi$ contains all primes (i.e. the property of being $\pi$\nobreakdash-isolated  coincides with the property of  being simply isolated), this statement was proved in the paper \cite{McDon}, and, in general case, in the work \cite{Bezv}. In fact, the assertion of Theorem stated above was only a part of the results in these works. Namely, in them it was proved  the existence of algorithm for constructing a system of generators for, respectivelly,  isolator and $\pi$-isolator of finitely generated subgroup of a free group and therefore  proofs in both papers were rather complicated and were carried out with using a number of nontrivial properties of the Nielsen set of generators of subgroups of free groups. Nevertheless, the proof of only the existence of a finite system of generators for isolators of  finitely generated subgroup of  free group proposed here is quite elementary and requires only three preliminary remarks.\medskip

{\bf Proposition 1.} {\it Let $H$ be a finitely generated subgroup of the free group $F$ and let the inclusion $f^{m}\in H$ hold for some nonidentity element $f\in F$ and some positive integer $m$. Then the rank of subgroup $K$ of group $F$ generated by subgroup $H$ and element $f$ does not exceed the rank of subgroup $H$.}\medskip

{\sl Proof.} Let $a_1, a_2, \dots, a_n$ be a system of free generators of subgroup $H$. By hypothesis, for some word $w = w(a_i)$ in these generators the equality $f^{m} = w$ holds. Note that since $f^{m}\neq 1$ the word $w$ is not empty.

Further, let $U$ be the free group with free generators $ x_1, x_2, \dots, x_n, y $ and let $V$ be its subgroup generated by elements $x_1, x_2, \dots, x_n$. Then the group $T =\langle x_1, x_2, \dots, x_n, y; \ w(x_i)=y^{m}\rangle$ is, on the one hand, the quotient group of the group $U$ by the normal closure of element $w (x_i)y^{-m}$, and, on the other hand, the free product of group $V$ and infinite  cyclic group generated by the element $y$, with the amalgamated subgroups generated by elements $w(x_i)$ and $y^{m}$.

Let $\varphi$ denote the natural homomorphism of group $U$ onto the quotient group $T$. Also, the mapping $x_i\mapsto a_i\, (i=1,2, \dots, n)$ and $y\mapsto f$ defines, obviously, a homomorphism $\psi$ of group $T$ onto the subgroup $K$ of group $F$. Then the product $\tau=\varphi\psi$ is an epimorphism of the group $U$ of rank $n + 1$ onto the free group $K$ and therefore the rank of $K$ does not exceed $n + 1$. But since the kernel of homomorphism $\tau$ is different from the identity subgroup and since free groups of finite ranks are Hopfian, the rank of the subgroup $K$ must be strictly less than $n + 1$. 

So, proposition 1 is proved.\medskip

The following statement is almost obvious:\medskip

{\bf Proposition 2.} {\it For any set $\pi$ of primes  and for arbitrary subgroup $H$ of some group $G$ any  subgroup $K$ of $\pi$-isolator $I_{\pi}(H)$ such that $H\leqslant K$ and $K\neq I_{\pi}(H)$ is not $\pi$-isolated.}\smallskip

Indeed, otherwise the inclusion $I_{\pi}(H)\subseteq K$ would hold.\medskip

The third statement we need was formulated in the book \cite {MagnusKarrasSolitar1974} (see p.114, Problem 2.4.17) as Takahasi's theorem: \medskip

{\bf Proposition 3.} {\it If $H_1\subseteq H_2\subseteq\dots$ is an ascending chain of subgroups each of rank $\leqslant r$ in a free group $F$ then all $H_i$ coincide for $i$ suficiently large.}\medskip

Going directly to the proof of the Theorem, suppose, arguing by contradiction, that the $\pi$-isolator $I_{\pi}(H)$ of some finitely generated subgroup $H$ of  free group $F$ is not a finitely generated subgroup. It will be shown that then in the group $F$ there exists an infinite strictly increasing sequence of subgroups such that the rank of each of them does not exceed the rank of the subgroup $H$. Thus, in view of Proposition 3 the Theorem will be proved.

Let us put $H_1=H $. Then subgroup $H_1$ is contained in the subgroup $I_{\pi}(H)$ and, being finitely generated, does not coincide with this subgroup and therefore is not $\pi$-isolated. So, there is an element $f_1$ that does not belong to the subgroup $H_1$ and is such that $f_1^{m_1}\in H_1$ for some positive $\pi$-number $m_1$. Since element $f_1^{m_1}$ belongs to the subgroup $I_{\pi}(H)$ (as it is $\pi$-isolated), we have the inclusion $f_1\in I_{\pi}(H)$. Therefore, the subgroup $H_2$ generated by subgroup $H_1$ and element $f_1$ is contained in $I_{\pi}(H)$, contains the subgroup $H_1$ but not coincides with it and is finitely generated. Moreover, it follows from Proposition 1 that its rank does not exceed the rank of the subgroup $H_1$.

Suppose now that for some $r\geqslant 2$ a strictly increasing sequence $H_1<H_2<\cdots<H_r$ of subgroups lying in $I_{\pi}(H)$ and such that the rank of each of them does not exceed the rank of the subgroup $H$ has already been constructed. Applying the arguments from the previous paragraph to the subgroup $H_r$, we find an element $f_r$ that does not belong to $H_r$ and is such that $f_r^{m_r}\in H_r$ for some positive $\pi$-number $m_r$. Then, as above, the subgroup $H_{r+1}$ generated by the subgroup $H_r$ and the element $f_r $ lies in  subgroup $I_{\pi}(H)$, contains the subgroup $H_r$  not coincides with it and has a rank not exceeding the rank of the subgroup $H$.

The existence of a sequence with the promised properties is proven.

\bigskip

 Ivanovo State University \smallskip

{\it E-mail address}: moldav@mail.ru


\begin{thebibliography}{99}

\bibitem{Bezv}
{\it Bezverkhnii V.\,N.} Î $\pi$-isolator of free group // Algoritmic problems in group and semigroup theory: Collection of scientific papers. Tula State Pedagogital Institute 1990. P.3--13.(Russian)

\bibitem{Kont}
{\it Kontorovich P.\,G.} Groups witth  basis of splitting II // Math. sb., 1946. V.19(61), ¹2, P.287-–308.(Russian)

\bibitem{MagnusKarrasSolitar1974}
 {\it Magnus W.,  Karrass A. and  Solitar D.} Combinatorial group theory. Interscience, New York, 1966.

\bibitem{McDon}
{\it McDonough T.\,P.} Root-Closure in Free Groups // J. London Math. Soc. 1970. P. 191-–192.

\end{thebibliography}
\end{document}